\documentclass[11pt, twoside, eqno]{article}
\usepackage{latexsym}
\usepackage{amssymb}
\usepackage{amsfonts}
\begin{document}
\begin{center}

\noindent {\bf \Large Comment ``On dual ordered semigroups" 
}\bigskip

\medskip

{\bf Niovi Kehayopulu}\bigskip

{\small Department of Mathematics,
University of Athens \\
15784 Panepistimiopolis, Athens, Greece \\
nkehayop@math.uoa.gr }
\end{center}
\date{ }

\noindent{\bf Abstract.} This is about the paper by Thawhat Changphas 
and Nawamin Phaipong in Quasigroups and Related Systems 22 (2014), 
193--200. {\small \bigskip

\medskip

}
\bigskip

\noindent This paper is actually the introduction and the main part 
from section 2 (: the first decomposition theorem) of the paper by 
{\it S. Schwarz} [On dual semigroups, Czech. Math. J. 10(85) (1960), 
201--230 ([10] in the References of the paper). The authors 
considered an ordered semigroup $(S,.,\le)$ instead of the semigroup 
$(S,.)$ considered by Schwarz, repeated the results by Schwarz 
exactly as they are given by Schwarz in [3] (instead of referring to 
[3] in the text), without any further explanation. The order plays a 
very very little role in Lemma 1.1(1), Lemma 1.1(4), Lemma 2.7 and 
Lemma 2.15(4) (the Lemma 2.8 and Corollaries 9 and 10 are immediate 
consequences of Lemma 2.7), but the proofs of Lemma 1.1(1) and Lemma 
1.1(4) are wrong. According to the proof of Lemma 1.1(1), the 
relation $y\le x$ implies $yA\subseteq xA$. Is it true? Why? There is 
no explanation in the paper and it is not true. Besides, in an 
ordered semigroup, $y\le x$ does not imply $yA\subseteq xA$ in 
general. Does $y\le x$ imply $My\subseteq Mx$ in the proof of Lemma 
1.1(4)? This is also wrong. In addition, the Lemma 1.1(4) does not 
need any proof since it is an immediate consequence of Lemma 1.1(1) 
and the corresponding result on semigroups by Schwarz. In Lemma 
1.1(4), the authors repeated the proof given by Schwarz. It might be 
noted that this proof can be drastically simplified without using any 
additional results (taking an element in $Sr(M)$, we can immediately 
show that this element belongs to $r(M)$). As far as the Lemma 
2.15(4) is concerned, if $(S,.)$ is a semigroup with $0$, $A$ an 
ideal of $S$ and $L$ a left ideal of $r(A)$, then $L$ is also a left 
ideal of $(S,.)$. This is due to Schwarz and the authors repeated the 
proof due to Schwarz and added that ``if $x\in L$ and $S\ni y\le x$, 
then $y\in r(A)$" while only that very little part (the proof of 
$y\in r(A)$) was needed to be added in Schwarz's proof to pass from 
the semigroup to the ordered semigroup (see the Lemma 2,2 in [3]). 
The aim of the present note is to indicate the mistakes and correct 
the proofs of Lemma 1.1(1) and Lemma 1.1(4). Throughout the paper the 
authors repeat the results due to Schwarz. We will mainly deal with 
the proofs related to the order, for the rest the authors should 
refer to [3]). If $(S,.\le)$ is an ordered semigroup, the zero of $S$ 
is an element of $S$ usually denoted by $0$, such that $x0=0x=0$ and 
$0\le x$ for every $x\in S$, that is, $0$ is the zero of the 
semigroup $(S,.)$ and the zero of the ordered set $(S,\le)$. If $S$ 
is a semigroup or an ordered semigroup with $0$, $l(A)$ denotes the 
subset of $S$ defined by $l(A)=\{x\in S \mid xA=0\}$ and $r(A)$ the 
subset of $S$ defined by $r(A)=\{x\in S \mid Ax=0\}$. It is more than 
clear that one can never say that ``Obtained results generalize (or 
extend) the results on semigroups without order" as it is mentioned 
in the abstract and in the introduction of the paper.\\
{\bf This is the Lemma 1.1(1) in [1]:} If $(S,.,\le)$ is an ordered 
semigroup and $A$ a nonempty subset of $S$, then the set $l(A)$ is a 
left ideal and the set $r(A)$ is a right ideal of $S$. (The authors 
should added that $0\in S$).\\
{\bf This is the Lemma 1.1(4) in [1]:} If $(S,.,\le)$ is an ordered 
semigroup and $M$ a right ideal  of $M$, then the set $r(M)$ is an 
ideal of $S$. (Again here the $0\in S$ should be added).\\
{\bf Concerning the proof of Lemma 1.1(1):} If $A\not=\emptyset$, 
$x\in l(A)$ and $S\ni y\le x$, we have to prove that $y\in l(A)$. 
According to the authors, $y\le x$ implies $yA\subseteq xA$ which is 
wrong.\\
{\bf Concerning the proof of Lemma 1.1(4):} If $M$ is a right ideal 
of $S$, $x\in r(M)$ and $S\ni y\le x$, we have to prove that $y\in 
r(M)$. According to the authors, $y\le x$ implies $My\subseteq Mx$ 
which is also wrong.\\
Let us prove that if $M$ is a right (resp. left) ideal of an ordered 
semigroup $S$, then $y\le x$ does not imply $My\subseteq Mx$ (resp. 
$yM\subseteq xM$) in general. This shows the mistake in Lemma 1.1(1) 
as well, as the right and the left ideals of an ordered semigroup $S$ 
are nonempty subsets of $S$. $\hfill\Box$\medskip

\noindent{\bf Example.} [2] Consider the ordered semigroup 
$S=\{a,b,c,d,f\}$ defined by the multiplication and the covering 
relation given below:\begin{center}
$\begin{array}{*{20}{c}}
.&\vline& a&\vline& b&\vline& c&\vline& d&\vline& f\\
\hline
a&\vline& a&\vline& b&\vline& a&\vline& a&\vline& a\\
\hline
b&\vline& a&\vline& b&\vline& a&\vline& a&\vline& a\\
\hline
c&\vline& a&\vline& b&\vline& c&\vline& a&\vline& a\\
\hline
d&\vline& a&\vline& b&\vline& a&\vline& a&\vline& d\\
\hline
f&\vline& a&\vline& b&\vline& a&\vline& a&\vline& f
\end{array}$
$$\prec=\{a,b),(c,a),(d,a)\}.$$\end{center}The set $M=\{a,c,d\}$ is a 
left ideal of $S$, $c\le a$ but $cM\nsubseteq aM$. The set 
$M=\{a,b,c,d\}$ is a right ideal of $S$, $c\le a$ but $Mc\nsubseteq 
Ma$. $\hfill\Box$\\
\noindent{\bf The corrected form of Lemma 1.1(1) and its proof:} {\it 
If $(S,.,\le)$ is an ordered semigroup with zero and $A$ a nonempty 
subset of $S$, then the set $l(A)$ is a left ideal and the set $r(A)$ 
is a right ideal of $S$.}\\
\noindent{\bf Proof.} Let $A$ be a nonempty subset of $S$. The set 
$l(A)$ is a left ideal of $(S,.)$ [1]. Let now $x\in l(A)$ and $S\ni 
y\le x$. Then $y\in l(A)$, that is, $yA=\{0\}$. Indeed: Let $z\in A$. 
Then $yz\le xz\in xA=\{0\}$, so $xz=0$, and $yz=0$. Since 
$yA\subseteq \{0\}$ and $yA\not=\{0\}$, we have $yA=\{0\}$. 
$\hfill\Box$\\
\noindent{\bf Remark.} If $A$ is a nonempty subset of $S$, then $y\le 
x$ implies $yA\subseteq (xA]$. Indeed: If $ya\in yA$ for some $a\in 
A$, then $ya\le xa\in xA$, so $ya\in (xA]$. If $x\in l(A)$ and $S\ni 
y\le x$, then  $y\in l(A)$, that is, $yA=\{0\}$. In fact: Since $y\in 
l(A)$ and $yA\subseteq (xA]$, we have $yA\subseteq (\{0\}]=\{0\}$, 
and so $yA=\{0\}$. $\hfill\Box$\\
\noindent{\bf The corrected proof of Lemma 1.1(4):} If $M$ is a right 
ideal of $(S,.,\le)$, then $r(M)$ is an ideal of $(S,.,\le)$. Indeed: 
By Lemma 1.1(1), the set $r(M)$ is a right ideal of $S$, the rest of 
the proof is a consequence of the Lemma 1,1 d) in [3]. 
$\hfill\Box$\\
\noindent{\bf This is the Lemma 2.15(4):} {\it Let S be an ordered 
semigroup having zero and $A$ an ideal of $S$. If L is a ideal of 
$r(A)$, then L is an ideal of S}. \\It is enough to prove that ``if 
$x\in L$ such that $y\le x$, they $y\in r(A)$". \\As only the fact 
that $y\in r(A)$ should be added in Schwarz' proof, the authors 
should emphasize it and write: ``if $x\in L$ such that $y\le x$ then, 
since $L\subseteq r(A)$, we have $y\in r(A)$". Several times, and in 
Lemma 2.15, the authors used the fact that $ar(A)=\{0\}$. This can be 
used in the proof of the Lemma 2.15(4) as well as follows: If $A$ a 
nonempty subset of $S$, $L\subseteq r(A)$, $x\in L$ and $S\ni y\le 
x$, then $Ay=0$, and so $y\in r(A)$. Besides, in the proof the 
``$\{0\}\cup L\subseteq L$" should be replaced by ``$\{0\}\cup L=L$", 
since $0\in L$. $\hfill\Box$\medskip

Let us discuss now some further results related to a semigroup 
$(S,.)$ in an attempt to show that the proof of Lemma 1,1 d) in [3] 
(and in [1]) can be simplified without using any additional results. 
We can directly prove that $Sr(M)\subseteq r(M)$ by taking an element 
of $Sr(M)$ and show that it belongs to $r(M)$ and so for the 
inclusion $A\subseteq r(l(A))$ in Lemma 1.1(2). \medskip

\noindent{\bf Remark.} This is the Lemma 1,1 d) in [3] (which is the 
same proof of Lemma 1.1(4) in [1]): If $(S,.)$ is a semigroup with 0 
and $M$ a right ideal of $S$, then $r(M)$ is an ideal of $S$. Its 
proof in [3] (and in [1]) is based on the following three steps:

(a) If $(S,.)$ is a semigroup with 0 and $M$ a right ideal of $S$, 
then $Mr(M)=\{0\}$.\\Indeed: Let $a\in M$ and $b\in r(M)$. Since 
$Mb=0$, we have $ab\in Mb=\{0\}$, then $ab=0$, so $Mr(M)\subseteq 
\{0\}$. Since $M$ is a right ideal of $S$, we have $0\in M$, and 
$\{0\}\subseteq Mr(M)$.

(b) If $(S,.)$ is a semigroup with $0$ and $M$ a right ideal of $S$, 
then $M(S(r(M))=\{0\}$ and $M(r(M)S)=\{0\}$. Indeed: By (a), we have 
$M(S(r(M))=(MS)r(M)\subseteq Mr(M)=\{0\}$, so $M(S(r(M))=\{0\}$ and 
$M(r(M)S)=(Mr(M))S=\{0\}S=\{0\}$.

(c) The equality $M(S(r(M))=\{0\}$ implies $Sr(M)\subseteq r(M)$ and 
$M(r(M)S)=\{0\}$ implies $r(M)S\subseteq r(M)$. Indeed: Suppose 
$M(S(r(M))=\{0\}$. If $ab\in Sr(M)$ for some $a\in S$, $b\in r(M)$, 
then $M(ab)\subseteq M(Sr(M))=\{0\}$, then $M(ab)\subseteq\{0\}$, 
$M(ab)=\{0\}$, $ab\in r(M)$. Suppose $M(r(M)S)=\{0\}$. If $ab\in 
r(M)S$ for some $a\in r(M)$, $b\in S$, then $M(ab)\subseteq 
M(r(M)S)=\{0\}$, then $M(ab)=\{0\}$, and $ab\in r(M)$.

As a consequence, if $(S,.)$ is a semigroup with 0 and $M$ a right 
ideal of $S$, then $r(M)$ (resp. $l(M)$) is an ideal of $M$. In a 
similar way we prove that if $M$ is a left ideal of $(S,.)$, then 
$l(M)$ is a left ideal of $S$.

A very simple proof of Lemma 1,1 d) in [3] is the following: \\
\noindent If $M$ a right ideal of $S$, then $Sr(M)\subseteq r(M)$. 
Indeed: Let $a\in S$ and $b\in r(M)$. Since $Mb=0$, we have 
$M(ab)=(Ma)b\subseteq (MS)b\subseteq Mb=0$, so $M(ab)\subseteq 
\{0\}$. Since $0\in M$, we get $\{0\}\subseteq M(ab)$, then 
$M(ab)=\{0\}$, and $ab\in r(M)$. $\hfill\Box$\medskip

\noindent{\bf Remark.} This is the Lemma 1.1(2): If $(S,.,\le)$ is an 
ordered semigroup, and $A$ a nonempty subset of $S$, then $A\subseteq 
r(l(A)$ and $A\subseteq l(r(A))$. This lemma is stated without proof 
in [3] (order plays no role in it), and this is its proof in [1]: 
Since $l(A)A=0$, we have $A\subseteq r(l(A))$. It is not necessary to 
have $l(A)A=0$, the equality $l(A)A=0$ has been used by Schwarz for 
other purpose. The proof is as follows: Let $x\in A$. Then $l(A)x=0$, 
and so  $x\in r(l(A))$. \medskip

\noindent Let us give an example of an ordered semigroup with $0$ in 
which the Lemma 1.1 of the paper in [1] can be applied: \medskip

\noindent{\bf Example.} [2] For the ordered semigroup $S$ defined by 
the multiplication and the covering relation below the element $a$ is 
the zero of $S$.\begin{center}
$\begin{array}{*{20}{c}}
.&\vline& a&\vline& b&\vline& c&\vline& d&\vline& f\\
\hline
a&\vline& a&\vline& a&\vline& a&\vline& a&\vline& a\\
\hline
b&\vline& a&\vline& a&\vline& a&\vline& b&\vline& c\\
\hline
c&\vline& a&\vline& b&\vline& c&\vline& a&\vline& a\\
\hline
d&\vline& a&\vline& a&\vline& a&\vline& d&\vline& f\\
\hline
f&\vline& a&\vline& d&\vline& f&\vline& a&\vline& a
\end{array}$
$$\prec=\{(a,b), (a,c),(a,d),(a,f)\}.$$\end{center}The right ideals 
of $S$ are the sets: $\{a\}$, $\{a,b,c\}$, $\{a,d,f\}$, $S$.\\The 
left ideals of $S$ are the sets: $\{a\}$, $\{a,b,d\}$, $\{a,c,f\}$, 
$S$.

\end{document}